\input amstex
\input amsppt.sty
\magnification=\magstep1
\hsize=30truecc
\baselineskip=16truept
\vsize=22.2truecm
\NoBlackBoxes
\nologo
\pageno=1
\topmatter
\TagsOnRight
\def\N{\Bbb N}
\def\Z{\Bbb Z}

\def\l{\left}
\def\r{\right}
\def\b{\bigg}

\def\({\b(}
\def\[{\b[}
\def\){\b)}
\def\]{\b]}

\def\t{\text}
\def\f{\frac}
\def\mo{\roman{mod}}
\def\em{\emptyset}

\def\eq{\equiv}

\def\ls{\leqslant}
\def\gs{\geqslant}
\def\al{\alpha}

\def\la{\lambda}
\def\Proof{\noindent{\it Proof}}
\def\Remark{\noindent{\it Remark}}

\hbox{Acta Arith. 129(2007), no.\,4, 397--402.}
\medskip
\title A characterization of covering equivalence\endtitle
\author Hao Pan (Shanghai) and Zhi-Wei Sun (Nanjing)\endauthor
\leftheadtext {H. Pan and Z.-W. Sun}

\abstract Let $A=\{a_s(\mo\ n_s)\}_{s=1}^k$ and $B=\{b_t(\mo\ m_t)\}_{t=1}^l$
be two systems of residue classes. If $|\{1\ls s\ls k:\,x\eq a_s\ (\mo\ n_s)\}|$
and $|\{1\ls t\ls l:\,x\eq b_t\ (\mo\ m_t)\}|$ are equal for all $x\in\Z$, then $A$ and $B$
are said to be covering equivalent.
In this paper we characterize the covering equivalence in a simple and new way.
Using the characterization we partially confirm a conjecture of R. L. Graham and K. O'Bryant.
\endabstract
\thanks 2000 {\it Mathematics Subject Classification}. Primary 11B25;
Secondary 11A07, 11B75.
\newline\indent
The second author is responsible for communications, and supported
by the National Science Fund (Grant No. 10425103) for
Distinguished Young Scholars in China.
\endthanks
\endtopmatter
\document

\heading {1. Introduction}\endheading

  For $n\in\Z^+=\{1,2,3,\ldots\}$ and $a\in\{0,\ldots,n-1\}$,
we simply use $a(n)$ to denote the residue class $\{x\in\Z:\, x\eq a\ (\mo\ n)\}$.
For a finite
system
$$A=\{a_s(n_s)\}_{s=1}^k\ \ (0\ls a_s<n_s)\tag1.1$$
of residue classes, the $n_1,\ldots,n_k$ are called its moduli and
its {\it covering function} $w_A:\Z\to\N=\{0,1,\ldots\}$ is given
by
$$w_A(x)=|\{1\ls s\ls k:\, x\in a_s(n_s)\}|.\tag1.2$$
(The covering function $w_{\em}$ of an empty system is regarded as
the zero function.) The periodic function $w_A(x)$ has many
surprising properties (cf. [S03a], [S04] and [S05a]).

Let $m$ be a positive integer. If $w_A(x)=m$ for all $x\in\Z$, then
(1.1) is said to be an {\it exact $m$-cover} of $\Z$ as in [S95] and [S96].
Recently Z. W. Sun (cf. [S04] and [S05b]) showed that $(1.1)$ forms an exact $m$-cover of $\Z$
if it covers $|S(n_1,\ldots,n_k)|$ consecutive integers exactly $m$ times, where
$$S(n_1,\ldots,n_k)=\l\{\f r{n_s}:\ r=0,\ldots,n_s-1;\, s=1,\ldots,k\r\}.\tag1.3$$
For problems and results on covers of $\Z$ by residue classes, the
reader is referred to [FFKPY], [G04] and [S03b].

 For two finite systems $A=\{a_s(n_s)\}_{s=1}^k$
and $B=\{b_t(m_t)\}_{t=1}^l$, Sun [S89] called $A$ and $B$ {\it covering equivalent}
(in short, $A\sim B$) if they have the same covering function (i.e., $w_A=w_B$).
Thus (1.1) is an exact $m$-cover of $\Z$ if and only if $(1.1)$ is covering equivalent to
the system consisting of $m$ copies of $0(1)$.

 In [S01] and [S02] Sun characterized the covering equivalence by various systems of equalities.
In this paper we present a simple characterization involving roots of unity.
Namely, we have the following result.

\proclaim{Theorem 1.1} Let $A=\{a_s(n_s)\}_{s=1}^k\ (0\ls a_s<n_s)$
and $B=\{b_t(m_t)\}_{t=1}^l$  $(0\ls b_t<m_t)$
be two finite systems of residue classes. Let $p$ be a prime greater than
$|S(n_1,\ldots,n_k,m_1,\ldots,m_l)|$,
and let $\zeta_p$ be a primitive $p$th root
of unity. Then $A$ and $B$ are covering equivalent
if and only if
$$\sum_{s=1}^k\f{\zeta_p^{a_s}}{1-\zeta_p^{n_s}}
=\sum_{t=1}^l\f{\zeta_p^{b_t}}{1-\zeta_p^{m_t}}.\tag1.4$$
\endproclaim
\proclaim{Corollary 1.1} $(1.1)$ forms an exact $m$-cover of $\Z$ if and only if
$$\sum_{s=1}^k\f{e^{2\pi i a_s/p}}{1-e^{2\pi in_s/p}}=\f m{1-e^{2\pi i/p}},\tag1.5$$
where $p$ is any fixed prime greater than $|S(n_1,\ldots,n_k)|$.
\endproclaim
\Proof. Simply apply Theorem 1.1 with $B$ consisting $m$ copies of $0(1)$. \qed
\medskip
\Remark\ 1.1. In 1975 \v S. Zn\'am [Z75a] used the transcendence of $e$ to prove that
$(1.1)$ is a disjoint cover (i.e., exact $1$-cover) of $\Z$ if and only if
$$\sum_{s=1}^k\f{e^{a_s}}{1-e^{n_s}}=\f1{1-e}.$$

\proclaim{Corollary 1.2} Suppose that for nonempty system $(1.1)$ we have
$$\sum_{s=1}^k\f{e^{2\pi ia_s/p}}{1-e^{2\pi in_s/p}}=0$$
where $p$ is a prime. Then
$$n_1+\cdots+n_k-k+1\gs|S(n_1,\ldots,n_k)|\gs p.\tag1.6$$
\endproclaim
\Proof. Clearly $|S(n_1,\ldots,n_k)|\ls n_1+\cdots+n_k-k+1$.
Since we don't have $A\sim\em$, applying Theorem 1.1 with $B=\em$ we find that
 $|S(n_1,\ldots,n_k)|$ cannot be smaller than $p$. This concludes the proof. \qed

 Corollary 1.2 partially confirms the following conjecture arising from the study of Fraenkel's
 conjecture on disjoint covers of $\N$ by Beatty sequences.
 \proclaim{Graham--O'Bryant Conjecture {\rm ([GO])}}
 Let $n_1,\ldots,n_k$ be distinct positive integers
 less than and relatively prime to $q\in\Z^+$. If $a_1,\ldots,a_k\in\Z$ and
 $$\sum_{s=1}^k\f{e^{2\pi ia_s/q}}{1-e^{2\pi in_s/q}}=0,$$
 then we must have $\sum_{s=1}^k n_s\gs q$.
 \endproclaim

 The following example shows that we cannot replace the prime $p$
 in Corollary 1.2 or Theorem 1.1
 by a composite number.

\medskip
 \noindent{\it Example} 1.1. Let $q>1$ be an integer and let $p$ be a prime divisor of $q$.
 Then, for any $n=1,\ldots,q-1$, we have
 $$\sum_{s=0}^{p-1}\f{e^{2\pi i(sq/p)/q}}{1-e^{2\pi in/q}}
 =\f{\sum_{s=0}^{p-1}e^{2\pi is/p}}{1-e^{2\pi in/q}}=0$$
 but $|S(n,\ldots,n)|=n<q$. Thus the conditions $0\ls a_s<n_s\ (s=1,\ldots,k)$ in Corollary 1.2
 cannot be cancelled.
 If $q$ is composite, then there are $q/p-1>0$ integers in the interval $((p-1)q/p,\, q-1]$.
 So we cannot substitute a composite number for the prime $p$ in Corollary 1.2.
\medskip

\proclaim{Corollary 1.3} Let $A=\{a_s(n_s)\}_{s=1}^k\ (0\ls
a_s<n_s)$ and $B=\{b_t(m_t)\}_{t=1}^l$  $(0\ls b_t<m_t)$ both have
distinct moduli. Let $p$ be a prime greater than
$|S(n_1,\ldots,n_k,m_1,\ldots,m_l)|$, and let $\zeta_p$ be a
primitive $p$th root of unity. Then $A$ and $B$ are identical if
and only if $(1.4)$ holds.
\endproclaim
\Proof. By a result of Zn\'am [S75b], $A$ and $B$ are identical if they have the same covering
function. Combining this with Theorem 1.1 we immediately get
the desired result. \qed

Observe that $A=\{a_s(n_s)\}_{s=1}^k$ and $B=\{b_t(m_t)\}_{t=1}^l$
are covering equivalent if and only if
$$\sum^k\Sb s=1\\x\in a_s(n_s)\endSb 1+\sum^l\Sb t=1\\x\in b_t(m_t)\endSb
(-1)=0\quad\ \t{for every}\ x\in\Z.$$
Thus Theorem 1.1 has the following equivalent form which will be proved in the next section.
\proclaim{Theorem 1.2} Let $\Cal A=\{\langle\la_s,a_s,n_s\rangle\}_{s=1}^k$ where
$\la_s,a_s, n_s\in\Z$ and $0\ls a_s<n_s$. Let $p>|S(n_1,\ldots,n_k)|$ be a prime,
and let $\zeta_p$ be any primitive $p$th root of unity. Then
$\Cal A\sim\em\ ($i.e., $w_{\Cal A}(x)=\sum_{1\ls s\ls k,\, x\in a_s(n_s)}\la_s=0$
for all $x\in\Z)$ if and only if
$$\sum_{s=1}^k\la_s\f{\zeta_p^{a_s}}{1-\zeta_p^{n_s}}=0.\tag1.7$$
\endproclaim

\heading{2. Proof of Theorem 1.2}\endheading

Let $S=S(n_1,\ldots,n_k)$. As $p>|S|\gs\max\{n_1,\ldots,n_k\}$,
there is a common multiple $N\in\Z^+$ of the moduli
$n_1,\ldots,n_k$ such that $N\eq 1\ (\mo\ p)$. Just as in [S05a],
we have
 $$\align\sum_{r=0}^{N-1}w_{\Cal A}(r)z^r
 =&\sum_{r=0}^{N-1}\sum\Sb 1\ls s\ls k\\n_s\mid a_s-r\endSb \la_sz^r
=\sum_{s=1}^k\la_s\sum\Sb 0\ls r<N\\r\in a_s(n_s)\endSb z^r
\\=&\sum_{s=1}^k\la_sz^{a_s}\sum_{0\ls q<N/n_s}(z^{n_s})^q
\\=&N\sum\Sb 1\ls s\ls k\\z^{n_s}=1\endSb\f{\la_s}{n_s}z^{a_s}
+(1-z^N)\sum\Sb 1\ls s\ls k\\z^{n_s}\not=1\endSb\la_s\f{z^{a_s}}{1-z^{n_s}}.
\endalign$$
Thus
$$\sum_{r=0}^{N-1}w_{\Cal A}(r)\zeta_p^r
=(1-\zeta_p^N)\sum_{s=1}^k\lambda_s\f{\zeta_p^{a_s}}{1-\zeta_p^{n_s}}.$$
It follows that
$$\sum_{s=1}^k\lambda_s\f{\zeta_p^{a_s}}{1-\zeta_p^{n_s}}=0
\iff\sum_{l=0}^{p-1}c_l\zeta_p^l=0,\tag2.1$$
where
$$c_l=\sum\Sb x=0\\x\in l(p)\endSb^{N-1}w_{\Cal A}(x)\in\Z.$$

 If $w_{\Cal A}(x)=0$ for all $x\in\Z$, then (1.7) holds by the above.

 Below we assume (1.7).
Then $\sum_{l=0}^{p-1}c_l\zeta_p^l=\sum_{r=0}^{N-1}w_{\Cal A}(r)\zeta_p^r=0$.
In the case $N=1$, it follows that $w_{\Cal A}(x)=w_{\Cal A}(0)=0$ for all $x\in\Z$.

Now suppose $N>1$. Clearly $N>p$ as $N\eq1\ (\mo\ p)$.
Since $1+x+\cdots+x^{p-1}=(x^p-1)/(x-1)$ is
the minimal polynomial of $\zeta_p$ over the field of rational numbers,
we must have
$c_0=c_1=\cdots=c_{p-1}$. (See also M. Newman [N71].)
Observe that if $x\in\Z$ then
$$w_{\Cal A}(x)=\sum_{s=1}^k\f{\lambda_s}{n_s}\sum_{r=0}^{n_s-1}e^{2\pi i\f{a_s-x}{n_s}r}
=\sum_{\al\in S}e^{-2\pi i\al x}\sum^k\Sb s=1\\\al n_s\in\Z\endSb
\f{\la_s}{n_s}e^{2\pi i\al a_s}.\tag2.2$$
(This trick appeared in [S91] and [S04].)
Since $|S|<p$, for each $l=0,\ldots,|S|$ we have
$$\aligned c_l=&\sum\Sb x=0\\x\in l(p)\endSb^{N-1}w_{\Cal A}(x)
=\sum_{\alpha\in S}\sum\Sb s=1\\\alpha n_s\in\Z\endSb^k\f{\lambda_s}{n_s}e^{2\pi i\alpha a_s}
\sum\Sb x=0\\x\in l(p)\endSb^{N-1}e^{-2\pi i\alpha x}
\\=&\sum_{\alpha\in S}e^{-2\pi i\alpha l}
\sum\Sb s=1\\\alpha n_s\in\Z\endSb^k\f{\lambda_s}{n_s}e^{2\pi i\alpha a_s}
\sum_{j=0}^{\l\lfloor (N-1-l)/p\r\rfloor}e^{-2\pi i\alpha pj},
\endaligned$$
where $\lfloor\cdot \rfloor$ is the greatest integer function.
If $l\in\{1,\ldots,|S|\}$ then
$$\l\lfloor\f{N-1-l}p\r\rfloor=\f{N-1}p+\l\lfloor\f{-l}p\r\rfloor=\f{N-1}p-1;$$
if $\alpha\in S\setminus\{0\}$ then
$$C(\alpha):=\sum_{j=0}^{(N-1)/p-1}e^{-2\pi i\alpha pj}
=\f{1-(e^{-2\pi i\alpha p})^{(N-1)/p}}{1-e^{-2\pi i\alpha p}}
=\f{1-e^{2\pi i\alpha}}{1-e^{-2\pi i\alpha p}}\not=0.$$
Let $c=c_0=\cdots=c_{p-1}$. By the above,
$$\sum_{\alpha\in S}e^{-2\pi i\alpha j}f(\alpha)=c$$
for every $j=0,\ldots,|S|-1$, where
$$f(0)=\f{N-1}p\sum_{s=1}^k\f{\lambda_s}{n_s}$$
and
$$f(\alpha)=e^{-2\pi i\alpha}C(\alpha)
\sum\Sb s=1\\\alpha n_s\in\Z\endSb^k\f{\lambda_s}{n_s}e^{2\pi i\alpha a_s}
\quad \ \t{for}\ \alpha\in S\setminus\{0\}.$$

 Let $\alpha_0=0,\,\al_1,\ldots,\alpha_{|S|-1}$ be all the distinct elements of $S$.
Now that
$$\sum_{t=0}^{|S|-1}e^{-2\pi i\alpha_tj}f(\alpha_t)=c\quad\ \t{for each}\ j=0,\ldots,|S|-1,$$
by Cramer's rule
$D_t=Df(\alpha_t)$ vanishes for every $t=1,\ldots,|S|-1$,
where $D=\det((e^{-2\pi i\al_t})^j)_{0\ls j,t<|S|}$
is of Vandermonde's type and hence nonzero.
Therefore
$$\sum_{\Sp s=1\\ \alpha n_s\in\Z\endSp}^k\f{\lambda_s}{n_s}e^{2\pi i\alpha a_s}=0
\quad\ \text{for all}\ \alpha\in S\setminus\{0\}$$
and hence
$w_{\Cal A}(x)=\sum_{s=1}^k\lambda_s/n_s$
for all $x\in\Z$ by (2.2). It follows that
$$0=\sum_{r=0}^{N-1}w_{\Cal A}(r)\zeta_p^r
=\sum_{s=1}^k\f{\lambda_s}{n_s}(1+\zeta_p+\zeta_p^2+\cdots+\zeta_p^{N-1})
=\sum_{s=1}^k\f{\lambda_s}{n_s}\cdot\f{1-\zeta_p^N}{1-\zeta_p}.$$
So $\sum_{s=1}^k\lambda_s/n_s=0$ and hence ${\Cal A}\sim\emptyset$.
We are done.

\bigskip

(Hao Pan) Department of Mathematics, Shanghai Jiaotong University, Shanghai 200240, P. R. China.
{\it E-mail address}: {\tt haopan79\@yahoo.com.cn}

\medskip

(Zhi-Wei Sun) Department of Mathematics, Nanjing University, Nanjing 210093, P. R. China.
{\it E-mail address}: {\tt zwsun\@nju.edu.cn}

\enddocument